\documentclass[11pt]{article}

\usepackage{amsfonts,amssymb,amsmath,amsthm,latexsym}
\usepackage{comment}
\usepackage{enumerate}
\usepackage{graphics,graphicx}
\usepackage{ulem} 
\usepackage{color} 

\theoremstyle{plain}

\newtheorem{thm}{Theorem}[section]
\newtheorem{prop}[thm]{Proposition}
\newtheorem{lem}[thm]{Lemma}

\newtheorem{conjecture}[thm]{Conjecture}
\newtheorem{dfn}[thm]{Definition}
\newtheorem{prob}[thm]{Problem}

\newcommand{\Hom}{\mathrm{Hom}}

\newcommand{\pr}{\mathbf{P}}

\newcommand{\xx}{\mathbf{w}}

\newcommand{\TT}{\mathcal{T}}
\newcommand{\ff}{\mathcal{F}}
\newcommand{\HH}{H}

\newcommand{\xo}{\xx_{_{\! 1}}}
\newcommand{\xt}{\xx_{_{\! 2}}}
\newcommand{\ps}{p_{_{\! S}}}
\newcommand{\pst}{p_{_{\! {S_2}}}}
\newcommand{\pso}{p_{_{\! {S_1}}}}

\newcommand{\pt}{p_{_{\! T}}}
\newcommand{\yyt}{\xx_{_{\! T}}}
\newcommand{\xxt}{\xx_{_{\! T}}}
\newcommand{\xxs}{\xx_{_{\! S}}}

\newcommand{\sq}{\mbox{\raisebox{0.2ex}{\tiny{$\, \square \,$}}}}
\newcommand{\bb}{\Big}
\newcommand{\old}[1]{}
\newcommand{\pf}{\noindent {\bf Proof. }}
\newcommand{\coh}{(\cdot
	\,  | \, h)}
\newcommand{\ga}{\alpha}
\newcommand{\gbk}{\beta_{_k}}
\newcommand{\uu}[1]{u_{_{\! #1}}}
\newcommand{\vv}[1]{v_{_{\! #1}}}
\newcommand{\bean}{\begin{eqnarray*}}
\newcommand{\eean}{\end{eqnarray*}}

\begin{document}

\author{
David Conlon\thanks{Mathematical Institute, Oxford OX2 6GG,
United Kingdom. Email: {\tt david.conlon@maths.ox.ac.uk}. Research
supported by a Royal Society University Research Fellowship.}
\and
Jeong Han Kim\thanks{School of Computational Sciences,  Korea Institute for Advanced Study (KIAS),  Seoul, South Korea. Email: {\tt jhkim@kias.re.kr}.~
The author was supported by a National Research Foundation of Korea (NRF) Grant funded by the Korean Government (MSIP) (NRF-2012R1A2A2A01018585) and KIAS Internal Research Fund CG046001. This work was partially carried out while the author was visiting Microsoft Research, Redmond and  Microsoft Research, New England.}
\and
Choongbum Lee\thanks{Department of Mathematics,
MIT, Cambridge, MA 02139-4307. Email: {\tt cb\_lee@math.mit.edu}. Research supported  
by NSF Grant DMS-1362326.}
\and
Joonkyung Lee\thanks{
Mathematical Institute, Oxford OX2 6GG, United Kingdom. 
Email: {\tt joonkyung.lee@maths.ox.ac.uk}. Supported by the ILJU Foundation of Education and Culture.
}
}

\title{Some advances on Sidorenko's conjecture}
\date{}

\maketitle

\begin{abstract}
A bipartite graph $H$ is said to have Sidorenko's property if the probability that the uniform random mapping from $V(H)$ to the vertex set of any graph $G$ is a homomorphism is at least the product over all edges in $H$ of the probability that the edge is mapped to an edge of $G$. In this paper, we provide three distinct families of bipartite graphs that have Sidorenko's property. First, using branching random walks, we develop an embedding algorithm which allows us to prove that bipartite graphs admitting a certain type of tree decomposition have Sidorenko's property. Second, we use the concept of locally dense graphs to prove that subdivisions of certain graphs, including cliques, have Sidorenko's property. Third, we prove that if $H$ has Sidorenko's property, then the Cartesian product of $H$ with an even cycle also has Sidorenko's property.
\end{abstract}

\section{Introduction}

For two graphs $H$ and $G$, a \emph{homomorphism} from $H$ to $G$
is a mapping $g: V(H)\rightarrow V(G)$ such that $\{g(v),g(w)\}$
is an edge in $G$ whenever $\{v,w\}$ is an edge in $H$. Let $\Hom(H,G)$
denote the set of all homomorphisms from $H$ to $G$. A beautiful
conjecture of Sidorenko~\cite{Sidorenko} (also made by Erd\H{o}s
and Simonovits~\cite{Simonovits} in a slightly different form) asserts
the following.
\begin{conjecture}\label{conj:sido}
For all bipartite graphs $H$ and all graphs $G$,
\begin{equation}
|\Hom(H,G)|\ge|V(G)|^{|V(H)|}\left(\frac{2|E(G)|}{|V(G)|^{2}}\right)^{|E(H)|}.\label{eq:sido}
\end{equation}
\end{conjecture}

That is, the probability that a uniform random mapping from $V(H)$ to $V(G)$ is a homomorphism is at least the product over all edges of the probability that the edge is mapped to an edge of $G$.

We say that a graph $H$ has Sidorenko's property if \eqref{eq:sido} holds for all graphs $G$. While Sidorenko \cite{Sidorenko} himself noted that the conjecture holds for some simple graphs such as trees, even cycles and complete bipartite graphs, a spate of recent work~\cite{CoFoSu, Hatami, KiLeLe, LiSz}, some of which we will describe below, has greatly expanded the class of graphs known to have Sidorenko's property. In this paper, we further this progress, providing three more families of graphs that satisfy the conjecture.

Our first family consists of those graphs which are amenable to a certain embedding 
process using branching random walks. Though the method can be applied more 
broadly, we will focus on using the technique to show that graphs admitting a particular type of tree decomposition have Sidorenko's property.
A \emph{tree decomposition} of a graph $H$, a concept introduced by Halin~\cite{Halin}
and developed by Robertson and Seymour~\cite{RobertSeymour}, 
is a pair $(\mathcal{F}, \TT)$ consisting of a family $\mathcal{F}$  of vertex subsets of $H$ and  a tree $\TT$ on 
vertex set $\mathcal{F}$ satisfying
\begin{enumerate}
\item $\bigcup_{X\in\mathcal{F}}X=V(H)$, 
\item for each $\{v,w\} \in E(H)$, there exists a set $X \in \mathcal{F}$ such that
$v, w \in X$, and
\item for $X,Y,Z\in \mathcal{F}$, $X\cap Y\subseteq Z$ 
whenever $Z$ lies on the path from $X$ to $Y$ in $\TT$.
\end{enumerate}
A \emph{strong tree decomposition} of a graph $H$ is a tree decomposition $(\mathcal{F},\TT)$ of $H$ 
satisfying the following two extra conditions:
\begin{enumerate}
\item 
The induced subgraphs $H[X]$, $X \in \ff$, are edge-disjoint trees.
\item 
For every pair $X,Y\in \ff$ which are adjacent in $\TT$, there is an isomorphism between the minimum subtrees of $H[X]$ and $H[Y]$ containing $X\cap Y$ that fixes $X \cap Y$. 
\end{enumerate}
A graph $H$ is \emph{strongly tree-decomposable} if it admits a strong tree decomposition.
One can easily verify that every strongly tree-decomposable graph is bipartite.
Our first theorem is as follows.

\begin{figure}
    \centering
    \includegraphics[width=0.4\textwidth]{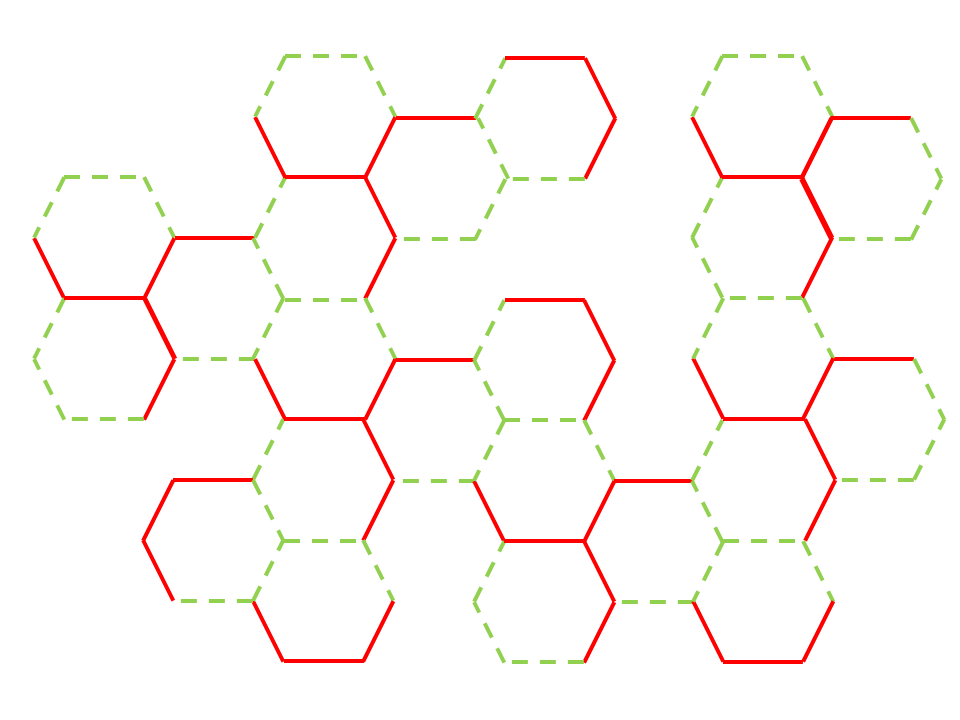}
    \caption{ A strongly tree-decomposable graph (each colored tree represents $H[X]$ for some $X\in\ff$).}\label{fig:hexagon}
\end{figure}
\begin{thm} \label{thm:main_thm}
If $H$ is strongly tree-decomposable, then $H$ has Sidorenko's property.
\end{thm}

The family of strongly tree-decomposable graphs completely contains some other families for which the conjecture was previously known.
For example, the family of reflection trees studied by Li and Szegedy~\cite{LiSz} consists exactly of those strongly tree-decomposable graphs
for which $\TT$ is a star, while the family of tree-arrangeable graphs considered in~\cite{KiLeLe} consists of those graphs $H$ such that $\mathcal{F}=\{N_H(a)\cup\{a\} \,|\, a\in A\}$ for an independent set $A$ in $V(H)$. However, our new family is broader than either of these, Figure \ref{fig:hexagon} showing an example of a strongly tree-decomposable graph that is neither a reflection tree nor tree arrangeable.

The second family of graphs we consider consists of subdivisions of certain graphs. The {\em subdivision}
of a graph $H$ is the graph obtained from $H$ by replacing its edges with paths of length two that are internally vertex disjoint. 
As a corollary of a more general theorem relating Sidorenko's conjecture to a conjecture about subgraphs of locally dense graphs, we obtain the following results.

\begin{thm} \label{thm:main_thm_2}
(i) For every positive integer $k$, the subdivision of $K_k$ has Sidorenko's property.

(ii) If $H$ has Sidorenko's property, then the subdivision of $H$ also has Sidorenko's property.
\end{thm}

The third family of graphs are those formed through taking certain Cartesian products.
The {\em Cartesian~product} $H_1 \sq H_2$ of two graphs $H_1$ and $H_2$ is the
graph on $V(H_1) \times V(H_2)$ where $(x_1, x_2)$
and $(y_1, y_2)$ are adjacent if and only if either (i) $x_1$ and $y_1$ are adjacent in $H_1$ and $x_2 = y_2$, or (ii) $x_2$ and $y_2$ are adjacent in $H_2$ and $x_1 = y_1$. In \cite{KiLeLe},
it was proved that if $H$ is a graph having Sidorenko's property and $T$ is a tree, then
$H \sq T$ also has Sidorenko's property. In particular, this result
implies that grid graphs of all dimensions have Sidorenko's property. 
Here, we prove that a similar result holds when trees are replaced by even cycles.

\begin{thm} \label{thm:main_thm_3}
Let $k \ge 2$ be a positive integer. If $H$ has Sidorenko's property, then $H \sq C_{2k}$ also has Sidorenko's property.
\end{thm}

In recent, independent work, Szegedy~\cite{Szegedy} developed a recursive procedure that generates a large family of graphs having Sidorenko's property. Szegedy's result is related to, but more general than, our Theorem~\ref{thm:main_thm}. For example, it shows that $2$-dimensional grids, which are not strongly tree-decomposable (but are amenable to other methods~\cite{KiLeLe}), have Sidorenko's property. Nevertheless, there are generalizations of Theorem~\ref{thm:main_thm}, obtained through further iteration analogous to the way in which one builds strongly tree-decomposable graphs from trees, which bring our result in line with Szegedy's. Since the statements of these more general results are quite technical, we have chosen to omit further discussion here. However, we refer the reader to the supplementary note~\cite{CoKiLeLe} for further information.

It is less clear how Szegedy's method relates to Theorems~\ref{thm:main_thm_2} and~\ref{thm:main_thm_3}. Szegedy has himself shown that his class is closed both under taking Cartesian products with trees and under taking subdivisions and it is plausible that it is also closed under taking Cartesian products with even cycles. This would give a satisfactory analogue of Theorem~\ref{thm:main_thm_3}, but, even if this doesn't hold, our result shows that one can simply add this new operation to the definition of the class. 

Parczyk~\cite{Parczyk} has also used Szegedy's method to recover Theorem~\ref{thm:main_thm_2}(i), that the subdivision of the complete graph has Sidorenko's property, but it is unlikely that his method can be used to recover all the results of Section~\ref{sec:subdiv} below, since they relate Sidorenko's conjecture to another conjecture whose status is undetermined. Things are muddied further by the fact that the results of Section~\ref{sec:subdiv} may also be generalized to hypergraphs. For example, if $H$ is a $k$-uniform hypergraph for which the analogue of Sidorenko's conjecture (or the Kohayakawa--Nagle--R\"odl--Schacht conjecture) holds, then we may subdivide each edge by adding a new vertex connected to each of the $k$ vertices in that edge to produce another graph which satisfies Sidorenko's conjecture. These considerations suggest the following open problem, which lies outside the reach of current methods.

\begin{prob}
Show that the subdivision of an arbitrary graph $H$ has Sidorenko's property.
\end{prob}

\section{Strongly tree-decomposable graphs} \label{sec:std}

The proof of Theorem \ref{thm:main_thm} is based on a   randomized 
algorithm that produces a homomorphism in $\Hom(H,G)$, with the desired bound
for $|\Hom(H,G)|$ following from an appropriate estimate on the entropy of the random
homomorphism.  

\medskip{}

Before getting into the proof, we give a brief overview of entropy.
Suppose that $X$ is a random variable taking values in some finite set $S$. 
The {\em entropy} of $X$ is defined by 
\[
	\HH(X)=\sum_{x\in S}\pr[X=x]\log\left(\frac{1}{\pr[X=x]}\right).
\]
Given two random variables
$X$ and $Y$ over the same space, 
the {\em relative entropy} of $Y$ with respect to $X$ is  
\[
	\HH(Y|X)=\HH(X,Y)-\HH(X).
\]
For three random variables $X,Y$, and $Z$, we say that $X$ and $Y$ are {\em conditionally independent}
given $Z$ if the distribution of $X$ conditioned on $(Y,Z) = (y,z)$ is identical to the distribution of $X$ conditioned on $Z=z$ for all possible values $y,z$. Finally, we note some simple properties of entropy which will be crucial in what follows.

\begin{lem}\label{lem:entropy}
Let $X$, $Y$, and $Z$ be random variables and suppose that $X$ takes values in 
a set $S$. Then 
\begin{enumerate}[(i)]
	\item $\HH(X)\leq\log|S|$,
	\item $\HH(X|Y,Z)=\HH(X|Z)$ if $X$ and $Y$ are conditionally independent given $Z$.
\end{enumerate}
\end{lem}

\subsection{Trees} \label{sec:treeEmbed}

We start with the case where $H = T$ is a tree. Let $r$ be a fixed vertex of $T$ and think of $T$ as a tree rooted at $r$. To generate a (not necessarily uniform) random homomorphism $\xx_T$ from $T$  to $G$, we consider the \emph{$T$-branching random walk}, or simply $T$-BRW,  on $G$. If we fix an order $r= u_{_{\! 0}}, u_{_{\! 1}}, .. , u_{_{\! t}}$ of all vertices in  $T$ in which no child precedes its parent, this algorithm proceeds as follows:
\begin{enumerate}
\item For the root $r$,  $\xxt (r)$ is a random vertex in $G$ chosen according to the stationary distribution, that is, with probability proportional to its degree,
$$
\pr[\xxt (r)=y]=\frac{\deg_G(y)}{2|E(G)|}.
$$
\item 
Suppose that $u_{_{\! i}}$, $i\geq 1$, is the first vertex that has not yet been embedded. Then, if $\vv{i}$ is the parent of $\uu{i}$, $\xx_T (u_{_i})$ is a vertex  chosen uniformly at random from all vertices that are adjacent to $\xxt (\vv{i})$ in $G$. (Note that $\vv i=\uu j$ for some $j<i$, as $\vv i$ must precede $\uu i$.) 

In others words, when the child $\uu{i}$ of $\vv{i}$ is born, it resides on a  vertex of $G$ chosen uniformly at random from all vertices that are adjacent to $\xxt(\vv{i})$ in $G$.   

\item Repeat Step 2 until all vertices have been embedded.
\end{enumerate}

Unless $G$ has no edge, $\xx_T$ is clearly a homomorphism from $T$ to $G$. Moreover, the distribution of $\xx_T$ does not depend on either the choice of the root vertex $r$ or the order of the remaining vertices. To see this, fix $h\in\Hom(T,G)$ and note that
\begin{align}\label{eq:root_invar_tree}
\pr[\xxt=h]&=
\frac{\deg_G(h(r))}{2|E(G)|}
\left(\frac{1}{\deg_G(h(r))}\right)^{d_r}
\prod_{v\in V(T)\setminus\{r\}}\left(\frac{1}{\deg_G(h(v))}\right)^{d_v-1}\nonumber \\ 
&=
\frac{1}{2|E(G)|}
\prod_{v\in V(T)}\left(\frac{1}{\deg_G(h(v))}\right)^{d_v-1},
\end{align}
where $d_v=\deg_T(v)$. 
As the last term depends only upon $T$ and $G$,  the claim follows.

We will call the distribution of $\xxt$ on $\Hom(T,G)$ the \emph{$T$-BRW distribution} and write $\pt:\Hom(T,G)\rightarrow [0,1]$ for the probability distribution  function of the $T$-BRW, i.e., 
$$\pt(h)=\pr[\xxt=h]. $$ 
Notice that (\ref{eq:root_invar_tree}) implies 
\begin{equation}\label{nz} 
\pt(h) >0 ~~\mbox{for all}~ h \in \Hom(T, G) .
\end{equation}
\old{We also note that the $T$-BRW for $T=K_2$
is a map that yields  a uniform random edge of $G$ and hence its distribution is the uniform distribution on $\Hom(K_2,G)$. }

The $T$-BRW also has an interesting projection property: if $S$ is a fixed subtree of $T$, then, for an arbitrary homomorphism $h : V(T) \rightarrow V(G)$, the map $h$ restricted to $V(S)$, which is denoted by $h\vert_{V(S)}$, or simply $h\vert_{S}$, induces a homomorphism from $S$ to $G$. Hence, $\xxt \vert_{S}$ induces another probability distribution over $\Hom(S,G)$. The following proposition asserts that $\xxt \vert_{S}$ is identically distributed to $\xxs$. For a subset $X $ of $V(T)$ and $h: X\rightarrow V(G)$, let { $\mathcal{G}_h (T) $} be the set of all homomorphisms from $T$ to $G$ that extend $h$, i.e., $\mathcal{G}_h (T) = \{ g\in\Hom(T,G):  g =h~{\rm on} ~ X\}$.

\begin{prop}\label{prop:tree_proj}
Let $S$ be a subtree of
 a tree $T$ and $h\in\Hom(S,G)$. 
 Then,
$$ \pr  \bb[ \xxt = h ~{\rm on}~ S \bb] = \ps (h),$$
where $\xxt$ is  the $T$-BRW on $G$, or equivalently
$$
\sum_{g\in\mathcal{G}_h(T)}\pt (g)=\ps (h).
$$ 
\end{prop}

\pf
As the distribution of $\xxt$ is independent of the root as well as the order of the vertices, one may regard $\xxt$ as a $T$-BRW that maps vertices in $S$ before any other vertices in $T$. Then $\xxt |_S $ and $\xx_S$ are constructed in the same way and the desired equality follows. 
\qed 

\bigskip

For ordered edges $(u,v)$  and $(s,t)$ of $T$ and $G$, respectively, one may apply  Proposition \ref{prop:tree_proj} to the tree $S$ consisting of the single edge $(u, v)$ and the homomorphism $h$ on $S$ with $(h(u),h(v))=(s, t)$ to  obtain 
$$ \pr \bb[  (\xxt (u), \xxt(v)) = (s,t) \bb]  = \ps (h) = \frac{\deg_G(s)}{2|E(G)|} \frac{1}{ \deg_G(s)} = \frac{1}{2|E(G)|}. $$
In other words, the  $T$-BRW $\xxt $ restricted to any edge $(u,v)$ of $T$ has the uniform distribution on the set of all ordered edges in $G$. Therefore,
\begin{equation}\label{unif}
\HH (\xxt (u), \xxt(v) ) = \log (2|E(G)|) . 
\end{equation}

We conclude this subsection by giving a lower bound for the entropy of the $T$-BRW $\xxt$, which in turn implies that trees satisfy Sidorenko's property. 

\begin{prop}\label{treecase} For a tree $T$ and a graph $G$, the $T$-{\rm BRW}   $\xxt$ on $G$
satisfies 
$$ \HH( \xxt) \geq \log \bb(|V(G)|^{|V(T)|} \left( \frac{2|E(G)|}{|V(G)|^2} \right)^{|E(T)|}\bb)$$
and, therefore, 
$$ |\Hom(T, G) | \geq  |V(G)|^{|V(T)|} \left( \frac{2|E(G)|}{|V(G)|^2} \right)^{|E(T)|}. $$
\end{prop}

\pf The second part easily follows from the first part and Lemma~\ref{lem:entropy}(i), which implies that
$$ \HH(\xxt) \leq \log |\Hom(T, G) | . $$
For the first part, observe, by the definition of conditional entropy and Lemma~\ref{lem:entropy}(ii), that  
\bean
\HH( \xxt ) 
&=&
\HH(\xxt(r)) +\sum_{i=1}^{t} \HH( \xxt(\uu i) |\xxt(\uu 0, ..., \uu{i-1}) )\\
&=&
\HH(\xxt(r)) +\sum_{i=1}^{t} \HH( \xxt(\uu i) |\xxt(\vv i )),
\eean
 where $\vv i$ is the parent of $\uu i$. Denote by $C_v$ the set of all children of $v$.
 Rearranging the last sum, we have that
 \bean
\HH( \xxt ) 
&=&
\HH(\xxt(r)) +\sum_{v\in V(T)}\sum_{u\in C_v} \HH( \xxt(u ) |\xxt(v  )) \\
&=&
\HH(\xxt(r)) +\sum_{v\in V(T)} \sum_{u\in C_v} \left(\HH( \xxt(u ), \xxt(v  ))-\HH (\xxt(v))\right).
\eean  
Since, by (\ref{unif}),  
$$ \sum_{v\in V(T)}\sum_{u\in C_v} \HH( \xxt(u ), \xxt(v  )) = \sum_{uv \in E(T)} \HH( \xxt(u ), \xxt(v  )) = |E(T)| \log (2|E(G)|)$$
and, by Lemma~\ref{lem:entropy}(i),
$$ 
-\HH(\xxt(r)) + \sum_{v\in V(T)}\sum_{u\in C_v} \HH (\xxt(v))
= \sum_{v\in V(T)} (\deg_T (v) -1) \HH (\xxt(v)) 
\leq (2|E(T)| - |V(T)|) \log |V(G)|, $$
we have that 
$$ 
\HH( \xxt ) 
\geq \log \bb( |V(G)|^{|V(T)|} \bb(\frac{2|E(G)|}{|V(G)|^2}\bb)^{|E(T)|} \bb) , 
$$
as desired. 
\qed 

\subsection{Branching random walks starting from a partial embedding}
\label{subsec:partialEmbed}

In order to embed strongly tree-decomposable graphs, we will first need to understand $T$-BRWs starting from a partial embedding. Let $T$ be a tree, $G$ a graph with at least one edge, and $X$ a subset of $V(T)$. Let $S$ be the minimum subtree of $T$ that contains $X$. If $X = \emptyset$, then  $S$ consists of a single fixed vertex of $T$. Take an ordering of the vertices of $T$ such that the vertices of $S$ precede all other vertices and, if the tree $T$ is seen as being rooted at the first vertex, no child precedes its parent. The {\em  $T$-branching random walk} on $G$, or simply $T$-BRW,  starting from a partial embedding $h : X \rightarrow V(G)$ is then defined as follows:

\begin{enumerate}
\item     
If  $\mathcal{G}_h (S) := \{g \in \Hom(S, G) : g|_{X} = h\}$ is empty, then there is no $T$-BRW starting with $h$. Otherwise, $\yyt(\cdot \,  | \, h):V(T) \rightarrow V(G)$  restricted to $S$ is defined  according to the distribution 
\begin{align*}
	\pr \Big[\, \yyt (\cdot
	\,  | \, h)  =y ~{\rm on} ~ S \, \Big]
	=
	\left\{
		\begin{array}{ll}	
		\frac{\ps(y)}{\sum_{g\in\mathcal{G}_h(S) }\ps(g)} 
		& \quad \text{if } y \in \mathcal{G}_h (S) \\
		0 & \quad \text{if }y \notin \mathcal{G}_h(S)~ .
		\end{array}
	\right.
\end{align*}
(Recall that $\ps$ is the probability distribution function of the $S$-BRW on $\Hom(S,G)$ and $\ps (g) > 0$ for all $g\in \Hom(S,G)$ as seen in (\ref{nz}) above). 
\item 
Once $\yyt(\cdot \,  | \, h)$ is defined on $S$, the embedding $\yyt(\cdot \,  | \, h)$ may be defined in the same way that was used for the construction of the $T$-BRW. That is, for the first vertex $u$ that has not yet been embedded and its parent $v$,  $\yyt(u \,  | \, h)$ is a vertex chosen uniformly at random from all vertices adjacent to $\yyt(v \,  | \, h)$ in $G$.  

\item Repeat Step 2 until all vertices have been embedded.
\end{enumerate}

If $X=\emptyset$, the $T$-BRW starting from $h: X \rightarrow V(G)$ is the usual $T$-BRW. Moreover, it is not difficult  to show that the $T$-BRW starting from $h$ has the same distribution as the usual $T$-BRW conditioned on $\xxt\vert_{X} = h$. We conclude this subsection by proving this fact.

\begin{prop} With the same notation as above and $x\in \Hom(T,G)$,
$$ 
\pr \Big[\, \yyt \coh
  =x ~{\rm on} ~ S \, \Big]= 
  \pr \bb[ \xxt = x ~{\rm on}~ S ~ \bb| \, \xxt = h ~{\rm on} 
~ X\bb]
$$
and 
$$
\pr   [\yyt\coh = x]
=  
\pr \bb[ \yyt =x ~\bb| ~
 \yyt  =h ~{\rm on} ~ X~\bb]. 
$$
In particular, the distribution of $\xxt\coh$ is independent 
of the order of vertices in $T$. 
\end{prop} 

\pf As Proposition \ref{prop:tree_proj} gives
$$ 
\ps (x|_S) = \pr [ \xxt = x ~{\rm on}~ S]$$
and 
$$ \sum_{g\in\mathcal{G}_h(S) }\ps(g)
= \sum_{g\in\mathcal{G}_h(S) }\sum_{f\in \mathcal{G}_{g} (T)}\pt(f)
= \sum_{f \in \mathcal{G}_{h} (T)}\pt(f)
= \pr [ \xxt = h ~{\rm on} 
~ X],$$
 we have that 
$$ \pr \Big[\, \yyt\coh  
  =x ~{\rm on} ~ S \, \Big] =
\frac{ \ps (x|_S)}{\displaystyle{\sum_{g\in\mathcal{G}_h (S)}\ps (g)}}
=\frac{\pr [ \xxt = x ~{\rm on}~ S]}{\pr [ \xxt = h ~{\rm on} 
~ X]}= \pr \bb[ \xxt = x ~{\rm on}~ S ~ \bb| \, \xxt = h ~{\rm on} 
~ X\bb] . 
$$
Since  $\xxt$ and  $\xxt\coh$ are constructed in the same way once their values on $S$ 
have been determined, we also have 
\begin{eqnarray*}
\pr   [\yyt\coh = x]
&=& 
\pr \Big[  \yyt\coh  =x ~{\rm on} ~ S \Big] 
\pr \Big[ \yyt\coh =x  \, \Big| \,  \yyt\coh   =x ~{\rm on} ~ S  \Big] \\
&=& \pr \bb[ \xxt = x ~{\rm on}~ S ~ \bb| \, \xxt = h ~{\rm on} 
~ X\bb]
\pr \Big[ \yyt  =x  \, \Big| \,  \yyt   =x ~{\rm on} ~ S  \Big] \\
&=& \pr \bb[ \xxt = x \, \Big| \,  \yyt   =h ~{\rm on} ~ X \Big],
\end{eqnarray*}
as required.
\qed 

\subsection{Embedding strongly tree-decomposable graphs}
\label{subsec:embedding}

Recall that a graph $H$ is strongly tree-decomposable if it 
has a tree decomposition $(\mathcal{F},\mathcal{T})$ satisfying two extra conditions:
\begin{enumerate}
\item 
The induced subgraphs $H[X]$, $X \in \ff$, are edge-disjoint trees.
\item 
For every pair $X,Y\in \ff$ which are adjacent in $\TT$, there is an isomorphism that fixes $X \cap Y$ between the minimum subtrees of $H[X]$ and $H[Y]$ containing $X\cap Y$. 
\end{enumerate}

Suppose that $G$ has at least one edge. We will attempt to construct a homomorphism  $\xx$ of $H$ into $G$ by regarding  $\TT$ as a rooted tree with root $R$, for a fixed $R\in \ff$, and taking a fixed order $R=X_0, X_1, X_2, ..., X_t$ of all sets in $\ff$ in which no child precedes its parent. We then proceed as follows: 

\begin{enumerate}
\item 
Take $\xx_{_{\!\!  H[R]}}$, the $H[R]$-BRW, and  set $\xx= \xx_{_{H[R]}}$ on $R$. Mark $R$ as explored. 

\item 
Suppose   $\xx$ is defined on $Z_{i-1}:=\bigcup_{j=0}^{i-1} X_i$,  $i\geq 1$, and $X_i$ is the first set that has not yet been explored. Take $\xx_{_{H[X_i]}} \coh$, the $H[X_i]$-BRW starting from $h=\xx|_{X_i\cap Z_{i-1}}$, and set $\xx=\xx_{_{H[X_i]}} \coh$ on $X_i\setminus Z_{i-1}$, provided that there exists a $H[X_i]$-BRW starting from $h$. Mark $X_i$ as explored.

If there is no $H[X_i]$-BRW starting from $h$, then the process stops and there is no output. (We will presently show that this case cannot occur.) 

\item 
Repeat Step 2 until every set in $\mathcal{F}$ has been explored.
\end{enumerate}

We first show that there is always an output $\xx$ and it is a homomorphism from $H$ to $G$. 

\begin{lem}\label{lem:orderInvariance}
The  above process  always outputs a homomorphism from $H$ to $G$. 
Moreover, the distribution of the  homomorphism does not depend on  the choice of the order 
$X_1, ..., X_t$ of the sets in $\ff$ provided that $X_0$ is  a fixed set $R$. 
\end{lem}

\pf 
We will prove inductively that $\xx$ is a homomorphism from $H[Z_i]$ to $G$ for $i=0,1,...,t$. 
For $i=0$, this follows since $G$ has at least one edge and $H[R]$ is a tree. Hence, 
the $H[R]$-BRW is well-defined and yields a homomorphism from $H[R]$ to $G$.
For $i\geq 1$,  suppose $\xx|_{Z_{i-1}} $ is a homomorphism. 
Since  $X_i$ is a leaf of the subtree $\TT_i:=\TT[ \{X_0, X_1, X_2, ..., X_i \}]$ of $\TT$,
every non-trivial path in  $\TT_i$ containing $X_i$ must contain its parent $Y_i$ and
the third property of tree decompositions gives 
$$ X_i \cap  Z_{i-1} = X_i \cap Y_i . $$

Since $H$ is strongly tree-decomposable, there is an isomorphism that fixes $X_i \cap Y_i$ 
between the minimum subtrees of $H[X_i]$ and $H[Y_i]$ containing $X_i\cap Y_i$. Since 
$Y_i \subset Z_{i-1}$, the relevant subtree of $H[Y_i]$ has already been
embedded by $\xx$. In particular, there is a homomorphism, say $g$, from 
the  relevant subtree of $H[X_i]$ to $G$ with  $g = \xx $ on $X_i \cap Y_i$.  
Hence, by the results of Section~\ref{subsec:partialEmbed}, there exists an $H[X_i]$-BRW starting from $\xx|_{X_i \cap Y_i}$. 
This shows that $\xx$ is a well-defined mapping from $Z_{i}=Z_{i-1} \cup X_i$ to $V(G)$. 

To show that $\xx$ is a homomorphism on $Z_i$, suppose $\{u, v\}$ is  an edge
of $H$ with $u,v \in Z_{i}$. Since both $\xx |_{_{Z_{i-1}}}$ and $\xx |_{_{X_i}}$ are 
homomorphisms, it is enough to show that either  $u, v \in Z_{i-1} $  or $u, v \in X_{i} $, though not necessarily exclusively. 
If not, we may assume that $u \in X_i$ and $ v\in X_j\subset Z_{i-1}$ for some $j <i$. 
We now take a set $U\in \ff$ that contains both $u$ and $v$. Such a set exists by the second property of tree decompositions. 
If $U$ is a descendant of $X_i$ in $\TT$, then the path from $U$ to $X_j$ 
in $\TT$ must go through $X_i$. As $v \in U \cap X_j$, we have 
$v \in X_i$ and hence $u,v \in X_i$, contradicting our assumption. 
If $U$ is not a descendant of $X_i$ in $\TT$, then 
the path from $X_i$ to $U$ in $\TT$ must contain $Y_i$, the parent of $X_i$. 
Since $ u \in X_i\cap U$, we have  $u \in Y_i \subset Z_{i-1}$ and hence 
$u,v \in Z_{i-1}$, again contradicting our assumption.

Therefore, the process constructing $\xx$ continues until $\xx$ has been extended to all sets in $\ff$. 
Since $\bigcup_{X\in \ff} X = V(H)$ (by the first property of tree decompositions), $\xx$ is a homomorphism
from $H$ to $G$. Moreover, using
\begin{align*} \label{eq:prob_explicit}
\pr\bb[\xx=h \bb]
&=\pr[\xx =h \text{ on }X_0]
  \prod_{i=1}^{t} \pr\left[\xx = h \text{ on }X_i ~\Big\vert~ \xx=h \text{ on } Z_{i-1}  \right]
\end{align*}
and 
$$
	\pr\left[\xx = h \text{ on }X_i ~\Big\vert~ \xx=h \text{ on } 
	Z_{i-1} \right]
	=
	\pr\left[\xx = h \text{ on }X_i ~\Big\vert~ \xx=h \text{ on } X_i \cap Y_i \right],
	$$
\newcommand{\gs}{\sigma}
we easily see that if $C_Y$ is the set of children of $Y$ in $\TT$, then
\begin{eqnarray*}
\pr\bb[\xx=h \bb]
&=&\pr\bb[\xx =h \text{ on }X_0 \bb]  \, 
  \prod_{i=1}^{t} \pr\left[
  \xx = h \text{ on }X_i ~\Big\vert~ \xx=h \text{ on } X_i \cap Y_i
   \right] \\
   &=&\pr\bb[\xx =h \text{ on }X_0\bb]  \, 
  \prod_{Y\in \ff}\prod_{ X\in C_Y  } \pr\left[
  \xx = h \text{ on }X ~\Big\vert~ \xx=h \text{ on } X \cap Y
   \right].   
\end{eqnarray*}
Similarly, for any other order $R=X_0=X'_0, X'_{1}, X'_{2}, ..., X'_{t} $ in which 
no child precedes its parent,
\begin{eqnarray*}
\pr\bb[\xx'=h \bb]
   =\pr\bb[\xx' =h \text{ on }X'_0\bb]  \, 
  \prod_{Y\in \ff}\prod_{X\in C_Y } \pr\left[
  \xx' = h \text{ on }X ~\Big\vert~ \xx'=h \text{ on } X \cap Y
   \right].   
\end{eqnarray*}
On $X_0 = X'_0 = R$,  both $\xx$ and $\xx'$ are $H[R]$-BRWs 
and, hence, 
$$ \pr[\xx' =h \text{ on }X'_0]  = \pr[\xx =h \text{ on }X_0]. $$
Further, for all adjacent pairs $X,Y \in \mathcal{F}$, conditioned on
$\xx' = \xx  = h$ on $X \cap Y $, 
both $\xx|_{X} $ and $\xx'|_{X}$ are 
$H[X]$-BRWs starting from $h$ and 
 $$ \pr\left[
  \xx' = h \text{ on }X  ~\Big\vert~ \xx'=h \text{ on } X \cap Y
   \right]
   = \pr\left[
  \xx = h \text{ on }X ~\Big\vert~ \xx=h \text{ on } X \cap Y \right].
  $$
Thus, $\pr[\xx=h ]=\pr [\xx'=h ]$, which means that 
the distribution of $\xx$  is independent of the order. 
\qed 

\bigskip

In the following proposition, we show that the choice of root $R$ is also irrelevant.

\begin{prop} \label{prop:marginal}
The distribution of the homomorphism $\xx$ described above does not depend on the choice of root $R$. In particular, $\xx|_X$ is the $H[X]$-BRW
for each $X \in \mathcal{F}$, i.e.,  
\[
\pr[\xx  =h ~{\rm on}~  X]=p_{_{\! H[X]}}(h|_{X})
\]
for $h\in\Hom(H,G)$.
\end{prop}

\pf
 Since   $X$ in the second part may be chosen as the root of $\TT$, 
the equality easily follows from the first part.
To prove the first part, let $\xo$ and $\xt$ be the random homomorphisms 
constructed using the above procedure with roots $Y_1, Y_2 \in \mathcal{F}$, respectively.  
We will show that $\xo$ and $\xt$ have the same distribution
whenever $Y_1$ and $Y_2$ are adjacent in $\TT$. 
The general case follows easily from repeated application of this special case.

Since the order of sets in $\mathcal{F}$ other than the root does not 
affect the distributions of $\xo$ and $\xt$, we may embed $Y_2$ immediately after $Y_1$ when $Y_1$ is the root.
If, on the other hand, $Y_2$ is the root, then $Y_1$ is to be embedded immediately after $Y_2$. Order all remaining sets in $\mathcal{F}$ in the exact same order for both cases.
It therefore suffices to prove that, 
for each $h \in \Hom(H,G)$,
\begin{align} \label{eq:marginal_key}
	&\, \pr[\xo = h \text{ on }Y_1] \cdot \pr[\xo = h \text{ on }Y_2 ~|~ \xo = h \text{ on } Y_1 \cap Y_2 ] \nonumber \\
	 & =\, \pr[\xt = h \text{ on }Y_2] \cdot \pr[\xt = h \text{ on }Y_1 ~|~ \xt = h \text{ on } Y_1 \cap Y_2].
\end{align}
We will first show that if $S_1$ and $S_2$ are the minimum subtrees  of $H[Y_1]$ and $H[Y_2]$ containing 
$Y_1\cap Y_2$, respectively, then
$$
\pr[\xo = h \text{ on }Y_2 ~|~ \xo = h \text{ on } Y_1 \cap Y_2]
= \frac{\pr[\xt = h \text{ on }Y_2]}{\displaystyle{\sum_{g\in \mathcal{G}_h (S_2)} \pst (g)}} , $$
and, similarly, 
$$
\pr[\xt = h \text{ on }Y_1 ~|~ \xt = h \text{ on } Y_1 \cap Y_2 ]
= \frac{\pr[\xo = h \text{ on }Y_1  ]}{\displaystyle{\sum_{g\in \mathcal{G}_h(S_1)} \pso (g)}} , $$ 
where $\mathcal{G}_h (S_i)  = \{g \in \Hom(S_i, G) : g =h~{\rm on}~Y_1\cap Y_2  \}$. 
The first equation  follows from noting that
\begin{eqnarray*}
\pr\bb[\xo = h \text{ on }Y_2 ~\bb|~ \xo = h \text{ on } Y_1 \cap Y_2 \bb]
&=& \pr\bb[\xo = h \text{ on }S_2 ~\bb|~ \xo = h \text{ on } Y_1 \cap Y_2 \bb] \\
& &  ~~\times \pr\bb[\xo = h \text{ on }Y_2 ~\bb|~ \xo = h \text{ on } S_2 \bb] \\
&=& \frac{\pst (h|_{S_2})}{\displaystyle{\sum_{g\in \mathcal{G}_h (S_2)} \pst (g)}} 
 \pr \Big[\xo = h  \text{ on }Y_2 \, \bb|  \,  \xo = h \text{ on } S_2 \bb], 
\end{eqnarray*}

and 
\begin{eqnarray*}
\pr \bb[ \xt  = h \text{ on } Y_2 \, \bb] 
&=&  \pr\bb[ \xt = h \text{ on } S_2 \bb ] 
\pr \bb[ \xt  = h \text{ on }Y_2 \, \bb| \, \xt = h \text{ on } S_2 \bb] \\
&=&  \pst (h|_{S_2})
\pr \Big[\xo = h  \text{ on }Y_2 \, \bb|  \,  \xo = h \text{ on } S_2 \bb],
\end{eqnarray*}
where, in the last equality, we used the fact that $\xo$ and $\xt$ are constructed in the same way once 
their values on $S_2$ are determined to be $h$. 
The same argument also works for the second equation. 

The desired equality  (\ref{eq:marginal_key}) follows since 
$$ \sum_{g\in \mathcal{G}_h (S_1)} \pso (g)
=\sum_{g\in \mathcal{G}_h (S_2)} \pst (g), $$
which, in turn, easily follows from the fact that there is an isomorphism from $S_1$ to $S_2$ that fixes $Y_1 \cap Y_2 $. 
\qed 

\bigskip

We are now in a position to prove Theorem~\ref{thm:main_thm}.

\bigskip

\noindent
{\bf Proof of Theorem~\ref{thm:main_thm}.}
Let $H$ be a strongly tree-decomposable graph with tree decomposition $(\mathcal{F}, \mathcal{T})$ and let $G$ be a graph with at least one edge. 
Using 
$$ \HH(\xx ) = \HH(\xx(R)) + \sum_{i=1}^{t}  \HH(\xx(X_i)| \xx(Z_{i-1}) ) $$
and, from Lemma~\ref{lem:entropy}(ii), 
$$\HH(\xx(X_i)| \xx(Z_{i-1}) ) = \HH(\xx(X_i)| \xx(X_i \cap Z_{i-1}) )= \HH(\xx(X_i)| \xx(X_i \cap Y_i) ), $$
where $Y_i$ is the parent of $X_i$ in the tree $\TT$
rooted at $R$,
we have   that 
\bean
\HH(\xx ) = \HH(\xx(R)) + \sum_{i=1}^{t}  
\HH(\xx(X_i)| \xx(X_i \cap Y_{i}) )
&=&  \HH(\xx(R)) + \sum_{Y\in \ff}\sum_{X\in C_Y} 
\HH(\xx(X)| \xx(X\cap Y) ) , 
\eean
where $C_Y$ is the set of children of $Y$. As 
$ \HH(\xx(X)| \xx(X\cap Y))= \HH(\xx(X)) - \HH( \xx(X\cap Y)), $
we also have that 
\begin{align}
\nonumber
\HH(\xx ) 
&= 
\HH(\xx(R)) + \sum_{Y\in \ff}\sum_{X\in C_Y} 
\HH(\xx(X)) - \HH( \xx(X\cap Y)\nonumber)\\ \nonumber
&= \sum_{X\in \ff} \HH(\xx(X)) - 
\sum_{XY\in E(\TT)} \HH( \xx(X\cap Y)) \\ 
&\geq  \sum_{X\in \ff} \HH(\xx(X)) - 
\sum_{XY\in E(\TT)} |X\cap Y| \log |V(G)| , \label{eqn:ApplyUpperBound}
\end{align}
where the inequality follows from Lemma~\ref{lem:entropy}(i).
Proposition \ref{prop:marginal} implies that, for each $X\in \ff$,
the distribution of $\xx|_{_X} $
is the $H[X]$-BRW distribution and, therefore, by Proposition~\ref{treecase},
$$ 
 \HH(\xx(X)) \geq \log 
\bb(|V(G)|^{|X|} \bb( \frac{2|E(G)|}{|V(G)|^2} \bb)^{|E(H[X])|}\bb)
= |E(H[X])| \log \frac{2|E(G)|}{|V(G)|^2} + |X|\log |V(G)|. $$
Since strong tree decomposability implies that
$$\sum_{X\in \ff} |E(H[X])| = |E(H)|   ~~~\mbox{and} ~~~ 
\sum_{X\in \ff} |X| -\sum_{XY \in E(\TT)} |X\cap Y| = |V(H)|, $$   
it follows that 
\bean
H(\xx) 
&\geq&
 \sum_{X\in \ff} 
|E(H[X])| \log 
  \frac{2|E(G)|}{|V(G)|^2} 
+\bb(\sum_{X\in \ff} |X| -\sum_{XY \in E(\TT)} |X\cap Y| \bb) \log |V(G)| \\
&=&
|E(H)|  \log 
  \frac{2|E(G)|}{|V(G)|^2} 
+|V(H)|  \log |V(G)| \\
&=&
\log 
\bb(|V(G)|^{|V(H)|} \bb( \frac{2|E(G)|}{|V(G)|^2} \bb)^{|E(H)|}\bb)~. 
\eean 
Since 
$H(\xx) \leq \log |\Hom(H,G)|$, 
we finally have that   
$$ \log |\Hom(H,G)| \geq  \log 
\bb(|V(G)|^{|V(H)|} \bb( \frac{2|E(G)|}{|V(G)|^2} \bb)^{|E(H)|}\bb), $$
or 
$$ |\Hom(H,G)| \geq 
 |V(G)|^{|V(H)|} \bb( \frac{2|E(G)|}{|V(G)|^2} \bb)^{|E(H)|} , $$
as desired. 
\qed

\section{Subdivisions} \label{sec:subdiv}

We say that a graph $G$ is $(\rho, d)$-dense if $G$ has density at least $d$ on every $U \subseteq V(G)$ with $|U| \geq \rho |V(G)|$. A beautiful conjecture of Kohayakawa, Nagle, R\"odl, and Schacht~\cite{KoNaRoSc} states that for any graph $H$ and any $\gamma, d > 0$, there should exist $\rho > 0$ such that any $(\rho, d)$-dense graph $G$ on a sufficiently large number of vertices contains at least $(1 - \gamma) |V(G)|^{|V(H)|} d^{|E(H)|}$ labeled copies of $H$. 

If $H$ is a bipartite graph which satisfies Sidorenko's conjecture, it clearly satisfies this conjecture with $\rho = 1$. However, the conjecture is also known for some non-bipartite graphs, including complete graphs, complete multipartite graphs, the line graph of the cube~\cite{CoHaPeSc}, and odd cycles~\cite{Reiher}. Here we show that any graph for which this conjecture holds may be subdivided to produce a graph that satisfies Sidorenko's conjecture.

\begin{thm} \label{thm:sec3}
Suppose that $H$ is a graph for which there exists a function $f: (0,1] \rightarrow (0, 1]$ and a constant $a > 0$ such that any $(f(d), d)$-dense graph $G$ contains at least
\[a |V(G)|^{|V(H)|} d^{|E(H)|}\]
copies of $H$. Then the subdivision of $H$ has Sidorenko's property.
\end{thm}

We note that both parts of Theorem~\ref{thm:main_thm_2} follow easily from this theorem and the discussion above. It will therefore suffice to prove Theorem~\ref{thm:sec3}. For our proof, it will be convenient to have a version of the hypothesis that applies to weighted graphs as well as graphs.

\begin{lem} \label{lem:weighttograph}
Suppose that $H$ is a graph for which there exists a function $f: (0,1] \rightarrow (0, 1]$ and a constant $a > 0$ such that any $(f(d), d)$-dense graph $G$ contains at least
\[a |V(G)|^{|V(H)|} d^{|E(H)|}\]
copies of $H$. Then, for any edge-weighted graph $W: V(G) \times V(G) \rightarrow [0,1]$ such that $\sum_{u, u' \in U} W(u, u') \geq d \binom{|U|}{2} $ for all $U \subseteq V(G)$ with $|U| \geq f(d/2)|V(G)|$, 
\[\sum_{x_1, \dots, x_h \in V(G)} \prod_{(i, j) \in E(H)} W(x_i, x_j) \geq \frac{a}{2}  |V(G)|^{|V(H)|} d^{|E(H)|}.\]
\end{lem}

We omit the proof of this lemma, since it is a standard application of Azuma's inequality (see, for example, Corollary~9.7 in \cite{CG15}).

The following lemma, due to Fox (see~\cite{KiLeLe}), says that in order to prove that a graph $H$ has Sidorenko's property, it suffices to prove that all graphs whose maximum degree is bounded in terms of the average degree contain roughly the correct number of copies of $H$.

\begin{lem} \label{lem:regular}
Let $H$ be a bipartite graph. If there exists a constant $c$, depending only on $H$, such that
\[|\Hom(H,G)|\ge c |V(G)|^{|V(H)|}\left(\frac{2|E(G)|}{|V(G)|^{2}}\right)^{|E(H)|}\]
for all graphs $G$ with maximum degree at most $\frac{4|E(G)|}{|V(G)|}$, then $H$ has Sidorenko's property.
\end{lem}

In proving Theorem~\ref{thm:sec3}, we may therefore assume that the maximum degree of $G$ is at most twice the average degree $d$. We will now argue that any graph of this type contains a large subgraph where the minimum degree is at least  $d/4$.

\begin{lem} \label{lem:mindeg}
If $G$ is a graph on $n$ vertices with average degree $d$ and maximum degree at most $2d$, there exists a subgraph $G'$ with at least $2^{-6} n$ vertices such that the minimum degree of $G'$ is at least $d/4$.
\end{lem}

\pf
We will prove the result by removing one vertex at a time to form graphs $G_0 = G, G_1, \dots, G_t$, eventually arriving at the required graph. Since the maximum degree of $G$ is at most $2d$, we see that $d(G_i) \leq 2d$ for all $i$, where $d(G_i)$ is the average degree of $G_i$. Suppose now that we have formed the graph $G_i$, but it does not satisfy the required condition. That is, there is a vertex $v_{i+1}$ of degree less than $d/4$.
We remove this vertex to form $G_{i+1}$. Since we start from a graph of average degree at least $d$, it follows that $d(G_i) \ge d$ for all $i$. Note that
\begin{align*}
 d(G_{i+1}) & = \frac{2|E(G_{i+1})|}{|V(G_{i+1})|} = \frac{2|E(G_{i})| - 2d(v_{i+1})}{|V(G_{i})| - 1} \geq \frac{2|E(G_{i})| - d/2}{|V(G_{i})| - 1} \\
 & =  d(G_i) + \frac{d(G_i) - d/2}{|V(G_i)| - 1} \geq d(G_i) + \frac{d}{2(|V(G_i)| - 1)}.
\end{align*}
Therefore,
\[d(G_t) \geq d + \frac{d}{2} \left(\frac{1}{n-1} + \frac{1}{n-2} + \dots + \frac{1}{n-t} \right),\]
which contradicts our observation that $d(G_t) \leq 2d$ when $t = (1 - 2^{-6}) n$. Hence, for some $i < t$, $G_i$ satisfied the required condition, completing the proof.
\qed

\bigskip

The key observation in our proof is the next lemma, which says that if a graph has large minimum degree, then the average codegree, taken over all pairs of vertices, in any large set $U$ is itself large.

\begin{lem} \label{lem:codegree}
Suppose that $G$ is a graph with $n$ vertices and minimum degree $\delta$. Then, provided $\delta|U| \geq 3n$,
\[\sum_{\substack{\{u, u'\} \subset U \\ u \neq u'}} d(u, u') \geq \frac{\delta^2}{2n} \binom{|U|}{2},\]
where $d(u, u')$ is the codegree of $u$ and $u'$. 
\end{lem}

\pf
If we write $d_U(v)$ for the number of neighbors of a vertex $v$ in $U$, we see, by convexity, that
\begin{align*}
\sum_{\substack{\{u, u'\}  \subset U \\ u \neq u'}} d(u, u') & = \sum_{v \in V(G)} \binom{d_U(v)}{2} \geq n \binom{\sum_v d_U(v)/n}{2}\\
& = n \binom{\sum_{u \in U} d(u)/n}{2} \geq n \binom{\delta|U|/n}{2} \geq \frac{\delta^2}{2n} \binom{|U|}{2},
\end{align*}
as required.
\qed

\bigskip

We now have all the necessary ingredients for proving Theorem~\ref{thm:sec3}.

\bigskip

\noindent
{\bf Proof of Theorem~\ref{thm:sec3}.}
Let $H_1$ be the subdivision of $H$. Suppose that $G$ is a graph with maximum degree at most $\frac{4|E(G)|}{|V(G)|}$.
By Lemma~\ref{lem:regular}, it will be enough to prove that there exists a constant $c$, depending only on $H$, such that
\begin{align} \label{eq:toshow}
|\Hom(H_1,G)|\ge c |V(G)|^{|V(H)| + |E(H)|}\left(\frac{2|E(G)|}{|V(G)|^{2}}\right)^{2 |E(H)|}.
\end{align}
Note that we may assume $|V(G)|$ is sufficiently large, since it is easy to see that verifying inequality~\eqref{eq:toshow} for a graph $G$ is equivalent to verifying it for a balanced blow-up of $G$.

By Lemma~\ref{lem:mindeg}, $G$ has a subgraph $G'$ with at least $2^{-6} |V(G)|$ vertices and minimum degree at least $\frac{|E(G)|}{2|V(G)|}$.
Now define a weighted graph $W$ on $V(G')$ by $W(u, v) = d(u, v)/|V(G')|$. By Lemma~\ref{lem:codegree}, the density of $W$ on every set $U \subseteq V(G')$ with $|U| \geq \frac{6 |V(G)|}{|E(G)|} |V(G')|$ is at least 
\[ 
	\frac{1}{{|U| \choose 2}} \sum_{\substack{\{u,u'\} \subseteq U\\ u \neq u'}} W(u,u')
	\ge \frac{1}{{|U|\choose 2}} \sum_{\substack{\{u,u'\} \subseteq U\\ u \neq u'}} \frac{d(u,u')}{|V(G')|}
		\ge 
	\frac{|E(G)|^2}{8 |V(G)|^4}. 
\]
Taking $d = |E(G)|^2/8 |V(G)|^4$, we have that $6|V(G)|/|E(G)| \leq f(d/2)$ for $|V(G)|$ sufficiently large.
Therefore, by our assumption on $H$ and Lemma~\ref{lem:weighttograph},
\[\sum_{x_1, \dots, x_h \in V(G')} \prod_{(i, j) \in E(H)} W(x_i, x_j) \geq \frac{a}{2}  |V(G')|^{|V(H)|} \left(\frac{|E(G)|^2}{8 |V(G)|^4}\right)^{|E(H)|}.\]
Expanding out $W$, this is the same as
\begin{align*}
\sum_{x_1, \dots, x_h \in V(G')} \prod_{(i, j) \in E(H)} d(x_i, x_j) & \geq a  |V(G')|^{|V(H)| + |E(H)|} \left(\frac{|E(G)|}{4 |V(G)|^2}\right)^{2 |E(H)|}\\
& \geq  2^{-6|V(H)| - 12 |E(H)|} a |V(G)|^{|V(H)| + |E(H)|} \left(\frac{2 |E(G)|}{|V(G)|^2}\right)^{2 |E(H)|}.
\end{align*}
But the first expression counts homomorphisms from $H_1$ to $G'$, while the last expression is of the required type \eqref{eq:toshow}. This completes the proof.
\qed

\bigskip

As an application of Theorem~\ref{thm:sec3}, we may prove the following more general theorem. We define the \emph{$K_{2,t}$-replacement} of a graph $H$ to be the graph obtained from $H$ by replacing each of its edges with a copy of $K_{2,t}$, identifying the edge's endpoints with the two vertices on one side of the corresponding copy of $K_{2,t}$. Note that a $K_{2,1}$-replacement is equivalent to a subdivision.

\begin{thm} \label{thm:replace}
Suppose that $H$ is a graph for which there exists an increasing function $f: (0,1] \rightarrow (0, 1]$ and a constant $a > 0$ such that any $(f(d), d)$-dense graph $G$ contains at least
\[a |V(G)|^{|V(H)|} d^{|E(H)|}\]
copies of $H$. Then the $K_{2,t}$-replacement of $H$ has Sidorenko's property.
\end{thm}

\pf
For a graph $G$ and an ordered sequence of $h = |V(H)|$ vertices $x_1, \dots, x_h$, let $H_1(x_1, \dots, x_h)$ be the number of homomorphisms from the subdivision $H_1$ of $H$ to $G$ such that the $i$th vertex of $H$ is mapped to $x_i$. By Theorem~\ref{thm:sec3}, 
\[\sum_{x_1, \dots, x_h} H_1(x_1, \dots, x_h) \geq |V(G)|^{|V(H)| + |E(H)|} \left(\frac{2 |E(G)|}{|V(G)|^2}\right)^{2 |E(H)|}.\]
Therefore, by convexity,
\begin{align*}
\sum_{x_1, \dots, x_h} H_1(x_1, \dots, x_h)^t & \geq |V(G)|^{|V(H)|} \left(\frac{\sum_{x_1, \dots, x_h} H_1(x_1, \dots, x_h)}{|V(G)|^{|V(H)|}}\right)^t\\ 
& \geq |V(G)|^{|V(H)| + t|E(H)|} \left(\frac{2 |E(G)|}{|V(G)|^2}\right)^{2 t|E(H)|}.
\end{align*}
Since the first expression counts homomorphisms from the $K_{2,t}$-replacement of $H$ to $G$, this completes the proof.
\qed

\bigskip

It is also possible to mix replacements in certain ways. For example, suppose that $H$ is a triangle and consider
\[X = \mathbb{E}_{x,y,z} d(x,y)^r d(y,z)^s d(z,x)^t.\]
This is a normalised count for the number of homomorphisms to $G$ of the graph formed by replacing the edges of the triangle by a $K_{2,r}$, a $K_{2,s}$, and a $K_{2,t}$, respectively. Note, by symmetry, that $X^3$ can be expressed in the form
\[X^3= \mathbb{E}[d(x,y)^r d(y,z)^s d(z,x)^t] \mathbb{E}[d(x,y)^s d(y,z)^t d(z,x)^r] \mathbb{E}[d(x,y)^t d(y,z)^r d(z,x)^s].\]
Hence, by H\"older's inequality with exponent $3$,
\[X \geq \mathbb{E}[(d(x,y)d(y,z)d(z,x))^{(r+s+t)/3}] \geq \mathbb{E}[d(x,y)d(y,z)d(z,x)]^{(r+s+t)/3} ,\]
where the last inequality follows by convexity. Therefore, the graph obtained by replacing the edges of the triangle by a $K_{2,r}$, a $K_{2,s}$, and a $K_{2,t}$ has Sidorenko's property. This argument easily generalises to all complete graphs and some other edge-transitive graphs. 

\section{Cartesian products with even cycles} \label{sec:cart}

Let $H, K$, and $G$ be  graphs.
In \cite{KiLeLe}, Kim, Lee, and Lee studied homomorphisms
from $H \sq K$ to $G$ by relating them to homomorphisms from
$H$ to the following auxiliary graph constructed from $K$ and $G$.

\begin{dfn}
Given graphs $K$ and $G$, let $\psi_K(G)$ be the graph
with vertex set $\Hom(K,G)$ in which two vertices 
$h_1, h_2 \in \Hom(K,G)$ are adjacent if and only if
$h_1(v)$ and $h_2(v)$ are adjacent in $G$ for all $v \in V(K)$.
\end{dfn}

The following lemma was proved in \cite{KiLeLe}. 

\begin{lem}\label{lem:1-1corresp}
For all graphs $H$, $K$ and $G$,
there is a one-to-one mapping between
$\Hom(H \sq K, G)$ and $\Hom(H, \psi_K(G))$.  
In particular,
\[
	|\Hom(H \sq K, G)| = |\Hom(H, \psi_K(G))|.
\]
\end{lem}

We are now ready to prove Theorem~\ref{thm:main_thm_3}.

\bigskip

\noindent
{\bf Proof of Theorem~\ref{thm:main_thm_3}.}
Let $k \ge 2$ be a given positive integer. 
It suffices to prove the theorem
when $H$ has no isolated vertex, since deleting isolated vertices in $H$
only removes connected components of $H\square J$ each of which are isomorphic to $J$.

Let   $G$ be  a  graph with $n$ vertices and let
$$p:=\frac{2|E(G)|}{n^2} = \frac{|\Hom(K_{2}, G)|}{n^2}, ~~
\ga := \frac{|\Hom(C_{4}, G)|}{p^{4} n^{4}},
~~~\gbk  := \frac{|\Hom(C_{2k}, G)|}{p^{2k} n^{2k}},
$$
for $k\geq 2$. 
By Lemma \ref{lem:regular}, 
we may assume that $G$ has maximum degree at most $2pn$.
Since $C_4$ and $C_{2k}$ have Sidorenko's property, $\alpha$ and  $\gbk$ are both at least $1$.
Using that $|\Hom(C_{2k}, G)| = \sum_{i=1}^{n} \lambda_i^{2k}$
for the eigenvalues   $\lambda_i$ of the adjacency matrix of $G$,
we also have
\[
	\gbk p^{2k} n^{2k}  \le  \alpha\lambda^{2k-4} p^4 n^4  
	~~~~~~\mbox{for $\lambda = \max_i |\lambda_i|$ and $k\geq 2$.}
\] 
Since $\lambda$ is at most the maximum degree of $G$, we know that $\lambda \le 2pn$. Thus the above yields
\begin{align} \label{eq:cartesian_0}
	\gbk \le 2^{2k-4} \alpha . 
\end{align} 
Writing $v = |V(H )|$ and $e = |E(H )|$, we apply Lemma \ref{lem:1-1corresp} and then use the fact that $H$ has Sidorenko's property to obtain 
$$ |\Hom(H \sq C_{2k}, G)| = |\Hom(H, \psi_{C_{2k}}(G))| 
	\ge |V(\psi_{C_{2k}}(G))|^{v} \left(\frac{2|E( \psi_{C_{2k}}(G))|}{|V(\psi_{C_{2k}}(G))|^2}\right)^{e} . \nonumber 
	$$
	As $|V(\psi_{C_{2k}}(G))| = |\Hom(C_{2k}, G)|= \gbk p^{2k}n^{2k}$, 
	it follows that  
\begin{align} \label{eq:cartesian_1}
	|\Hom(H \sq C_{2k}, G)| 
	\geq  (\gbk p^{2k}n^{2k})^{v} \left(
	\frac{2|E( \psi_{C_{2k}}(G))|}{(\gbk p^{2k}n^{2k})^2}\right)^{e}. 
\end{align}
Lemma \ref{lem:1-1corresp} also yields 
$$ 2|E( \psi_{C_{2k}}(G))| = |\Hom(K_2, \psi_{C_{2k}}(G))| =
|\Hom(K_2 \sq C_{2k}, G)|
	= 
	|\Hom(C_{2k}, \psi_{K_2}(G))|.$$
Since  $C_{2k}$ has Sidorenko's property, we obtain  
\begin{align*}
	2|E( \psi_{C_{2k}}(G))|  
	&\ge
	|V(\psi_{K_2}(G))|^{2k} \left(\frac{2|E( \psi_{K_2}(G))|}{|V(\psi_{K_2}(G))|^2} \right)^{2k} \\
	&=	
 \bb(\frac{ 2|E( \psi_{K_2}(G))|}{|\Hom(K_2, G)|}\bb)^{2k}.  
\end{align*}
Applying Lemma \ref{lem:1-1corresp} once more, we have 
\[
2 |E( \psi_{K_2}(G))|=	|\Hom(K_2, \psi_{K_2}(G))| = 
|\Hom(K_2\sq K_2, G )| = |\Hom(C_4, G)| = \alpha p^4 n^4 ,
\]
which, together with $|\Hom(K_2, G)| = 2|E(G)| = pn^2$, yields 
\[
	2|E( \psi_{C_{2k}}(G))|   
	\geq
	\bb( \frac{\alpha  p^4 n^4 }{pn^2} \bb)^{2k}
	=
	\alpha^{2k} p^{6k} n^{4k} .
\]
Together with~\eqref{eq:cartesian_1}, this inequality implies that 
\[
	|\Hom(H \sq C_{2k}, G)|
	\ge
	(\gbk p^{2k} n^{2k} )^{v} \bb(\frac{\alpha^{2k} p^{2k}}{\gbk^2}\bb)^{e}
	=
	\alpha^{2ek} \gbk^{v-2e} p^{2kv + 2ke} n^{2kv} .
\]
Since $H$ has no isolated vertex, we have $v - 2e \le 0$
and hence, by~\eqref{eq:cartesian_0}, $\gbk^{v-2e} \geq 2^{(2k-4)(v-2e)} \ga^{v-2e}$. 
Therefore, since $\ga \geq 1$,  
\[
	|\Hom(H \sq C_{2k}, G)|
	\ge
	2^{(2k-4)(v-2e)} \alpha^{2ek + (v-2e)} n^{2kv} p^{2kv + 2ke}
	\geq 2^{(2k-4)(v-2e)}  n^{2kv} p^{2kv + 2ke}.
	\]
Since $H \sq C_{2k}$ has $2kv$ vertices and $2kv + 2ke$ edges, the result now follows from Lemma~\ref{lem:regular}.
\qed 

\vspace{3mm}
\noindent
{\bf Acknowledgements.}
We would like to thank Olaf Parczyk for sharing his Master's thesis with us. We would also like to thank Bal\'azs Szegedy for a number of valuable discussions.

\end{document}